\documentclass[11pt]{article}
\usepackage{latexsym,amssymb,amsmath,amsfonts,amsthm}
\usepackage{graphics}
\usepackage{enumitem}
\usepackage{makecell}
\usepackage{graphicx}
\usepackage{algorithmicx}
\usepackage{algorithm}
\usepackage{mathrsfs}
\usepackage{subfigure}
\usepackage{color}
\usepackage{hyperref}
\usepackage{cite}
\usepackage{accents}
\makeatletter
\def\widebar{\accentset{{\cc@style\underline{\mskip14mu}}}}
\makeatother
\hypersetup{hidelinks}
\newcommand{\R}{{\mathbb R}}

\newcommand{\C}{{\mathbb C}}

\newcommand{\be}{\begin{eqnarray}}
\newcommand{\ben}{\begin{eqnarray*}}
\newcommand{\en}{\end{eqnarray}}
\newcommand{\enn}{\end{eqnarray*}}

\newcommand{\ov}{\overline}

\newcommand{\g}{\gamma}

\newcommand{\wi}{\widehat}
\newcommand{\wid}{\widetilde}
\newtheorem{theorem}{Theorem}[section]

\newtheorem{remark}[theorem]{Remark}
\definecolor{lk}{rgb}{0,0,1}
\begin{document}

\renewcommand{\theequation}{\arabic{section}.\arabic{equation}}

\title{\bf A deep learning approach to inverse medium scattering: Learning regularizers
from a direct imaging method}
\author{Kai Li\thanks{Academy of Mathematics and Systems Science, Chinese Academy of Sciences,
Beijing 100190, China and School of Mathematical Sciences, University of Chinese Academy of Sciences,
Beijing 100049, China ({\tt likai98@amss.ac.cn})}
\and
Bo Zhang\thanks{State Key Laboratory of Mathematical Sciences and Academy of Mathematics and Systems Science,
Chinese Academy of Sciences, Beijing 100190, China and School of Mathematical Sciences, University of Chinese
Academy of Sciences, Beijing 100049, China ({\tt b.zhang@amt.ac.cn})}
\and
Haiwen Zhang\thanks{Corresponding author. State Key Laboratory of Mathematical Sciences and Academy of Mathematics and Systems Science, Chinese Academy of Sciences, Beijing 100190, China ({\tt zhanghaiwen@amss.ac.cn})}
}
\date{}

\maketitle

\begin{abstract}
This paper aims to solve numerically the two-dimensional inverse medium scattering problem with far-field data.
This is a challenging task due to the severe ill-posedness and strong nonlinearity of the inverse problem.
As already known, it is necessary but also difficult numerically to employ an appropriate regularization
strategy which effectively incorporates certain a priori information of the unknown scatterer to overcome
the severe ill-posedness of the inverse problem.
In this paper, we propose to use a deep learning approach to learn the a priori information of the support
of the unknown scatterer from a direct imaging method.
Based on the learned a priori information, we propose two inversion algorithms for solving the inverse problem.
In the first one, the learned a priori information is incorporated into the projected Landweber method.
In the second one, the learned a priori information is used to design the regularization functional for the
regularized variational formulation of the inverse problem which is then solved with a traditional iteration
algorithm. Extensive numerical experiments show that our inversion algorithms provide good reconstruction results
even for the high contrast case and have a satisfactory generalization ability.

\vspace{.2in}

\end{abstract}

\section{Introduction}
\setcounter{equation}{0}

This paper is concerned with the inverse problem of scattering of time-harmonic acoustic waves from an
inhomogeneous medium in two dimensions.
This type of problem arises in various applications, such as sonar detection, remote sensing,
geophysical exploration, medical imaging, and nondestructive testing (see, e.g., \cite{C19}).

It is well-known that inverse medium scattering problems are strongly nonlinear and severely ill-posed.
A wide range of numerical reconstruction algorithms (including iterative algorithms and
non-iterative algorithms) have been developed, employing various regularization strategies, to recover
the inhomogeneous medium (or the contrast of the inhomogeneous medium) from a knowledge of the far-field
data or scattered-field data. For example, a continuation method was proposed in \cite{BL05} to recover the
inhomogeneous medium from multi-frequency scattering data, by first using the Born approximation to obtain
an initial guess of the inhomogeneous medium from the data at the lowest frequency and then applying
the Landweber method recursively on multiple frequencies to refine the reconstructed result.
A preconditioning technique was introduced in \cite{L10} for the iteratively regularized Gauss-Newton method
(IRGNM) and applied to solve inverse medium scattering problems.
The contrast source inversion (CSI) method was introduced in \cite{V97} for inverse medium scattering problems.
The basic idea of CSI involves minimizing the cost functional by alternatively updating both the contrast
and the contrast source.
A subspace-based optimization method (SOM) was proposed in \cite{C09}, which shares several properties
with CSI. See the monographs \cite{Chen2018,C19} and the references quoted therein for a comprehensive account
of various inversion methods for inverse scattering problems.
It should be noted that these iterative algorithms may have difficulties in obtaining
satisfactory reconstruction results for inhomogeneous media of high contrast.
One of the reasons for this may be that these iteration methods are probably not equipped with appropriate
regularization strategies which are actually difficult to be chosen in practical applications \cite{Z20,AJO17}.
For a comprehensive discussion of various regularization methods for ill-posed and inverse problems,
see the monographs \cite{E96,KNS08} and the references quoted therein.

Recently, non-iterative algorithms have attracted much attention in inverse medium scattering problems
since they do not need to solve the forward problem and thus are computationally fast.
For example, the orthogonality sampling algorithm was developed in \cite{P10} for the detection of the location
and shape of the unknown scatterer from the far field pattern, and an approximate factorization method was
introduced in \cite{ZZ20} for locating the support of the contrast of the inhomogeneous medium.
Nevertheless, non-iterative algorithms mainly focus on recovering the location and shape of the unknown scatterer,
and their reconstruction results are usually less accurate than those of iterative algorithms.
This paper aims to integrate the features of iterative algorithms and non-iterative algorithms via
a deep learning approach, thereby improving their reconstruction results.

In recent years, deep learning has been applied to develop powerful methods to solve various inverse
problems with very promising performance, such as computed tomography (CT) \cite{E24,G18,J17,HLJ20},
magnetic resonance imaging (MRI) \cite{Yang16,Yang20}, optical diffraction tomography (ODT) \cite{Y20},
electrical impedance tomography (EIT) \cite{W19,WZC19} and inverse medium scattering problems
(see, e.g., \cite{K19,L18,S19,GRL21,LTN22,LZZ24,LZR22,LYZ22,NJH24,Meng2024,WC19,WC18,Sun18,Z20}).
The reader is referred to \cite{M17,AMOS19,H23,LZ24,L23,Liang20,C20} for a good survey of deep learning-based
methods for various inverse problems including MRI image reconstruction problems (linear inverse problems) and
inverse medium scattering problems (strongly nonlinear inverse problems).
It should be mentioned that deep learning has been used to address the ill-posedness issue of inverse problems
by learning the undetermined regularization functionals from data for the variational regularization formulation
of some linear inverse problems (see, e.g., \cite{Yang16,Yang20,G18,HLJ20,Y20}).
However, to the best of our knowledge, \cite{LZZ24} is so far the only work in this direction for nonlinear
inverse scattering problems, where a deep learning-based iterative reconstruction algorithm was proposed
to solve inverse medium scattering problems, based on a repeated application of a convolutional neural network
(CNN) and the IRGNM. In \cite{LZZ24} we were not able to learn the undetermined regularization functional from
data directly; instead, we reformulated the regularized variational formulation with an unknown regularization
functional of the inverse problem as an equivalent constrained minimization problem with an
unknown feasible region depending on the undetermined regularization functional.
The CNN in \cite{LZZ24} is then designed to learn the a priori information of the shape of the unknown contrast
by using a normalization technique in the training process and trained to act like a projector which is helpful
for projecting the solution into the feasible region of the constrained optimization problem associated with
the inverse problem; see \cite{LZZ24} for details.

In this paper, we propose two iterative regularization algorithms that incorporate the a priori information
of the shape and location (i.e., the support) of the unknown contrast, which is learned by a deep
neural network from a direct imaging method, as regularization strategies for recovering the contrast
of the inhomogeneous medium from the far-field data.
Specifically, we first train a deep neural network to retrieve the a priori information of the support of the unknown contrast from a direct imaging method (e.g., the one in \cite{P10}).
Then, the learned a priori information is incorporated into the projected Landweber method in our first algorithm,
whilst the learned a priori information is used to construct the regularization functional
for the variational regularization formulation of the inverse problem which is finally solved by an iteration
algorithm in our second algorithm.
It is worth noting that the trained deep neural network in this paper is used to provide a good
approximation of the support of the unknown contrast (which is indeed confirmed in the numerical examples),
while the trained deep neural network in our previous work \cite{LZZ24} focused on learning the a priori
information of the shape of the unknown contrast.
%Therefore the iteration algorithms proposed in this paper should be equipped with more appropriate
%a priori information for the unknown contrast than the algorithm in \cite{LZZ24} and thus can achieve
%better reconstruction results in the high contrast case.
Extensive numerical experiments demonstrate that our algorithms have a satisfactory reconstruction
performance, strong robustness to noise and good generalization ability.

The rest of this paper is organized as follows. Section \ref{S2} presents the direct and inverse medium
scattering problems considered in this paper. In Section \ref{S3},
we propose a deep neural network that can retrieve
the support of the unknown contrast based on the direct imaging method.
In Section \ref{S4}, we present two iterative reconstruction algorithms that incorporate the a priori
information of the support of the unknown contrast for solving the inverse medium scattering problem.
Numerical experiments are carried out in Section \ref{S5} to illustrate the effectiveness of our algorithms.
Some conclusions and remarks are given in Section \ref{S6}.

\section{Problem formulation}\label{S2}
\setcounter{equation}{0}

In this section, we introduce the direct and inverse medium scattering problems considered in this paper.
Consider an inhomogeneous medium in $\R^2$ characterized by the piecewise smooth refractive index $n(x)>0.$
Define $m(x):=n(x)-1$ which is the contrast of the inhomogeneous medium and assumed to be compactly
supported in a disk with radius $\rho$, i.e., $\mathrm{supp}(m)\subset B_\rho:=\{x\in\mathbb{R}^2:|x|<\rho\}.$
Let $u^i=u^i(x,d):=e^{ikx\cdot d}$ be an incident plane wave with the incident direction
$d\in\mathbb{S}^1:=\{x\in\mathbb{R}^2:|x|=1\}$ and the wave number $k>0$.
Then the total field $u=u^i+u^s$, which is the sum of the incident field $u^i$ and the scattered field $u^s$,
satisfies the reduced wave equation
\begin{equation}\label{1}
\triangle u(x) +k^2n(x)u(x)=0 \qquad \text{in}\;\mathbb{R}^2,
\end{equation}
and the scattered field $u^s$ is assumed to satisfy the Sommerfeld radiation condition
\begin{equation}\label{2}
\lim_{r\to \infty}r^{\frac{1}{2}}\left(\dfrac{\partial u^s}{\partial r} - iku^s\right) = 0, \qquad r = |x|.
\end{equation}
It has been shown that the direct scattering problem \eqref{1}--\eqref{2} is well-posed (see, e.g., \cite{C19}).
Moreover, the scattered field $u^s$ has the following asymptotic behavior \cite{C19}
\begin{equation}\label{3}
u^s(x)=\dfrac{e^{ik|x|}}{\sqrt{|x|}}\biggl\{u^\infty(\hat{x})+\mathcal{O}\bigg(\dfrac{1}{|x|}\bigg)\biggr\},
\qquad \mathrm{as}\ |x|\to\infty,
\end{equation}
uniformly for all directions $\hat{x}=x/|x|$ on $\mathbb{S}^1$, where $u^\infty$ is the far-field pattern of
the scattered field $u^s$.
In the rest of the paper, we write the far-field pattern, the scattered field and the total field as $u^\infty(x,d)$,
$u^s(x,d)$ and $u(x,d)$, respectively, to indicate their dependence on the incident direction $d\in\mathbb{S}^1$.
It has been proved in \cite{C19} that the scattering problem \eqref{1}--\eqref{2} is equivalent to the well-known
Lippmann-Schwinger equation
\begin{equation}\label{7}
u(x, d) = u^i(x,d) +  k^2\int_{\mathbb{R}^2}\Phi(x,y)m(y)u(y, d)dy,\quad x\in \mathbb{R}^2,
\end{equation}
where $\Phi(x,y) := (i/4)H^{(1)}_0(k|x-y|)$, $x\neq y$, is the fundamental solution to the Helmholtz equation
$\Delta w +k^2 w=0 $ in two dimensions and $H^{(1)}_0$ denotes the Hankel function of the first kind of order zero.
Since $\mathrm{supp}(m)\subset B_\rho$, it is known from \cite{C19} that for any $d\in\mathbb{S}^1$, the integral
equation \eqref{7} is uniquely solvable in $C(\ov{B_\rho})$ and thus
	\be\label{hw-eq1}
	u=(I-k^2 T_m)^{-1}u^i
	\en
with $u=u(\cdot,d)$ and $u^i=u^i(\cdot,d)$,
where the operator $T_m: C(\ov{B_\rho})\rightarrow C(\ov{B_\rho})$ is given by
	\ben
	(T_m w)(x):=\int_{B_\rho}\Phi(x,y)m(y)w(y)dy,\quad x\in\ov{B_\rho},
	\enn
and $I-k^2 T_m:C(\ov{B_\rho})\rightarrow C(\ov{B_\rho})$ is bijective and has a bounded inverse.
Here, the subscript in $T_m$ indicates the dependence on the contrast function $m$.
Further, it is known from \cite{C19} that the far-field pattern $u^\infty$ of the scattered field $u^s=u-u^i$
is given by
\begin{equation}\label{8}
u^\infty(\hat{x},d)=\dfrac{k^{\frac{3}{2}}e^{i\frac{\pi}{4}}}{\sqrt{8\pi}}\int_{B_\rho}
e^{-ik\hat{x}\cdot y}m(y)u(y,d)dy
\end{equation}
for all $ \hat{x} = x/|x| $ on the unit circle $\mathbb{S}^1$.
We refer to \cite{C19} for details on the study of the direct scattering problem \eqref{1}--\eqref{2}.
In this work, we focus on the following inverse problem.

\textbf{Inverse problem (IP)}: Given the measured data $u^\infty(\hat{x},d)$ for $\hat{x}$, $d\in\mathbb{S}^1$,
reconstruct the unknown contrast $m$.

To present the mathematical formulation of the problem (IP), we introduce the far-field operator
$\mathcal{F}:L^2(B_\rho)\to L^2(\mathbb{S}^1\times\mathbb{S}^1)$ which maps the contrast $m(x)$
to its corresponding far-field pattern $u^\infty(\hat{x},d)$:
\begin{equation}\label{4}
\mathcal{F}(m)=u^\infty.
\end{equation}
From \eqref{hw-eq1} and \eqref{8} it follows that for $\hat{x},d\in\mathbb{S}^1$,
\be\label{hw-eq6}
(\mathcal{F}(m))(\hat{x},d)=\dfrac{k^{\frac{3}{2}}e^{i\frac{\pi}{4}}}{\sqrt{8\pi}}
\int_{B_\rho}e^{-ik\hat{x}\cdot y}m(y)((I-k^2 T_m)^{-1}u^i)(y)dy
\en
with $u^i=u^i(\cdot,d)$. From the representation \eqref{hw-eq6} it is seen that $\mathcal{F}$ is strongly nonlinear.
Further, since $e^{-ik\hat{x}\cdot y}$ is analytic for $\hat{x}\in\mathbb{S}^1$, it is clear that
the equation \eqref{4} is severely ill-posed.
For uniqueness results of the inverse scattering problem considered in this paper, the reader is referred to \cite{C19}.
Since only noisy measured data can be acquired in practice, we will consider the noised perturbation $u^{\infty,\delta}$
of the far-field data $u^\infty$ in the sense that
$\|u^{\infty,\delta}-u^\infty\|_{L^2(\mathbb{S}^1\times\mathbb{S}^1)}
\le\delta\|u^\infty\|_{L^2(\mathbb{S}^1\times\mathbb{S}^1)}$, where $\delta>0$ is called the noise level.
Thus, given the noisy measured data $u^{\infty,\delta}$, \eqref{4} can be rewritten as
\begin{equation}\label{5}
\mathcal{F}(m) \approx u^{\infty,\delta}
\end{equation}
for the unknown contrast $m$.

\begin{remark} {\rm
In the case when the contrast function $m$ is sufficiently small, then $u^i(\cdot,d)$ can be viewed as
an approximation of $u(\cdot,d)$, due to \eqref{7} and the method of successive approximations (see, e.g., \cite{C19}).
Thus we obtain the well-known Born approximation
\begin{equation}\label{9}
u^\infty(\hat{x}, d) \approx (\mathcal{F}_b(m))(\hat{x},d):=\dfrac{k^{\frac{3}{2}}
e^{i\frac{\pi}{4}}}{\sqrt{8\pi }}\int_{B_\rho}e^{-ik\hat{x}\cdot y}m(y)u^i(y, d)dy.
\end{equation}
It is clear that $\mathcal{F}_b$ is a linear operator from $L^2(B_\rho)$ to $L^2(\mathbb{S}^1\times\mathbb{S}^1)$.
Hence for the inverse problem (IP) in this case, the contrast $m$ can be approximately computed by solving the
linear equation \eqref{9} (see, e.g., \cite[Section 11.1]{C19}).
However, in this paper, we consider the contrast of the general form, which includes the high contrast case.
Thus we cannot apply the Born approximation \eqref{9} for the inverse problem (IP).
}
\end{remark}

In order to reconstruct the unknown contrast $m$, we need to compute the Fr\'{e}chet derivative of
the far-field operator $\mathcal{F}$.
\cite[Theorem 11.6]{C19} characterized the Fr\'{e}chet derivative of $\mathcal{F}$ in the case of three dimensions.
Similarly, it can be shown that the far-field operator $\mathcal{F}: m\mapsto u^\infty$ considered in this paper
is also Fr\'{e}chet differentiable with the derivative being given by $\mathcal{F}'(m)(q)=v^\infty $,
where $q\in L^2(B_\rho)$ and $v^\infty:=v^\infty(\hat{x},d)$ is the far-field pattern of the scattered field
$v^s$ satisfying the Sommerfeld radiation condition \eqref{2} and the reduced wave equation
\begin{equation}\label{6}
\triangle v^s + k^2nv^s = -k^2uq\qquad \textrm{in}\;\;\mathbb{R}^2.
\end{equation}
Here, $u = u(x,d)$ is the total field corresponding to the contrast $ m $.
From the above argument, it is found that one needs to compute the numerical solution of the equation \eqref{6}
for numerically solving the Fr\'{e}chet derivative of $\mathcal{F}$.

For the sake of numerical reconstruction, it is necessary to discretize the contrast $m$.
Specifically, define $C_\rho:=[-\rho,\rho]\times[-\rho,\rho]\subset\mathbb{R}^2$ and we discretize
$C_\rho$ into uniformly distributed $(N\times N)$ pixels which are denoted as $x_{ij}$ ($i,j=1,\ldots,N$).
Then the contrast $m$ can be approximately represented by a piecewise constant, which can be denoted by a
discrete matrix $\boldsymbol{m}=(\boldsymbol{m}_{ij})\in\mathbb{C}^{N\times N}$ with
$\boldsymbol{m}_{ij}:=m(x_{ij})$ and is also called the contrast matrix in the rest of the paper.
We further introduce a discrete matrix $ \boldsymbol{S} = (\boldsymbol{S}_{ij})\in \mathbb{R}^{N\times N} $ with
\begin{equation}\label{26}
\boldsymbol{S}_{ij} :=
\left\{ \begin{array}{ll}
1, & \boldsymbol{m}_{ij}\neq 0,\\
0, & \boldsymbol{m}_{ij} = 0
\end{array} \right.
\end{equation}
to characterize $ \mathrm{supp}(m) $. Note that the index set $\{(i,j):\boldsymbol{S}_{ij} = 1\}$ can be
viewed as the discrete form of $\mathrm{supp}(m)$. We also call $\boldsymbol{S}$ the support matrix in this paper.
Suppose the inhomogeneous medium is illuminated by $Q$ incident plane waves $u^i(x,d_q)$ with distinct
incident directions $d_q$ ($q=1,\dots, Q$) uniformly distributed on $\mathbb{S}^1$
and the far-field pattern is measured at $P$ distinct observation directions $\hat{x}_p$ ($p=1,\dots,P$)
uniformly distributed on $\mathbb{S}^1$.
The noisy far-field pattern $u^{\infty,\delta}$ can then be discretized as a measurement matrix
$\boldsymbol{u}^{\infty,\delta}:=(a^\delta_{p,q})\in\mathbb{C}^{P\times Q}$ with
$a^\delta_{p,q} := u^{\infty,\delta}(\hat{x}_p, d_q),\; p=1,\dots, P,\; q = 1,\dots, Q$.
Note that $\textrm{supp}(m)\subset B_\rho\subset C_\rho $.
The formula (\ref{5}) can be approximated as follows
\begin{equation}\label{24}
\boldsymbol{F}(\boldsymbol{m}) \approx \boldsymbol{u}^{\infty,\delta},
\end{equation}
where $\boldsymbol{F}$ denotes the discrete form of the far-field operator $\mathcal{F}$.
Denote by $\boldsymbol{F}'$ the discrete form of $\mathcal{F}'$.

\begin{remark}\label{R2} {\rm
%\color{lk}
The Lippmann-Schwinger equation \eqref{7} can be numerically solved efficiently by applying the fast
Fourier transform within a disk containing the support of the contrast $m$, as suggested by Vainikko
(see \cite{V00,H01}). In this paper, $\boldsymbol{F}$ and $\boldsymbol{F}'$ are computed by
this method with the disk to be $B_\rho$.
}
\end{remark}

\section{Deep learning for retrieving the support information}\label{S3}
\setcounter{equation}{0}

In this section, we propose a deep neural network to approximate the support of the unknown contrast,
based on the characterization of the support of the unknown contrast retrieved
by a direct imaging method. The retrieved information of the approximate support of the unknown
contrast will then be used in the next section to develop two regularized iteration algorithms for
solving the inverse problem (IP).

\subsection{The direct imaging method}\label{S3_1}

The key idea of the direct imaging method is to construct an appropriate imaging function in terms
of the measurement data to characterize the shape and location (i.e., the support) of the unknown
scattering obstacles or unknown contrast of an inhomogeneous medium (see, e.g., \cite{P10,ZZ20}).
In this paper, we adopt the following imaging function proposed in \cite[equation (33)]{P10} to
retrieve the information of the support of the unknown contrast $m$ from the far-field data
$u^\infty(\hat{x},d)$ at a fixed wave number $k$:
\begin{equation}\label{12}
I(z;u^\infty):=\int_{\mathbb{S}^1}\bigg|\int_{\mathbb{S}^1}e^{ik\hat{x}\cdot z}u^\infty(\hat{x},d)ds(\hat{x})\bigg|^2ds(d), \quad z\in \mathbb{R}^2.
\end{equation}
As discussed in \cite{P10}, it is reasonable to expect that the imaging function $I(z; u^\infty)$
takes a large value as the imaging point $z$ approaches the boundary of the support $\mathrm{supp}(m)$
and decays as $z$ moves away from $\mathrm{supp}(m)$.

As mentioned in Section \ref{S2}, only the noisy far-field data $u^{\infty,\delta}$ can be obtained
in practice. In order to visualize the imaging function of the direct imaging method numerically,
we define
\begin{equation}\label{25}
I^A(z; \boldsymbol{u}^{\infty,\delta}):=\dfrac{2\pi}{Q}\cdot\dfrac{2\pi}{P}\sum_{q=1}^{Q}
\left|\sum_{p=1}^{P}u^{\infty,\delta}(\hat{x}_p,d_q)e^{ik\hat{x}_p\cdot z}\right|^2
\end{equation}
for the measurement matrix $\boldsymbol{u}^{\infty,\delta}$ of the unknown contrast $m$ at a fixed
wave number $k$, where $P,Q$ are the same as in Section \ref{S2}.
Then $I^A(z;\boldsymbol{u}^{\infty,\delta})$ is a good trapezoid quadrature approximation to
the continuous imaging function $I(z;u^{\infty,\delta}).$
We further discretize $C_\rho=[-\rho,\rho]\times[-\rho,\rho]$ into uniformly distributed
$N\times N$ pixels, denoted as $z_{ij},\;i,j=1,\ldots,N$.
Then $I^A$ can be approximately represented by a discrete matrix
$\boldsymbol{I}=(\boldsymbol{I}_{ij})\in\mathbb{R}^{N\times N}$ with
$\boldsymbol{I}_{ij}:=I^A(z_{ij};\boldsymbol{u}^{\infty,\delta})$ and is also called
the imaging matrix of $m$ in the remaining part of the paper.

\subsection{Deep neural network for retrieving the approximate support}\label{S3_2}

In this subsection, we propose a deep neural network $\mathcal{M}_\Theta$ to extract some a priori
information of the support $\mathrm{supp}(m)$ of the unknown contrast $m$ (i.e., the support matrix
$\boldsymbol{S}$ defined in Section \ref{S2}) from the direct imaging method (i.e., the imaging
matrix $\boldsymbol{I}$ defined in Subsection \ref{S3_1}).
Here, $\Theta$ denote the parameters of the network to be determined during the training process.
Once it is well trained, $\mathcal{M}_\Theta$ should provide a good approximate support matrix for the
unknown contrast, which serves as important a priori information for reconstruction (see Section \ref{S4}
for more details). In this subsection, we will introduce the architecture of the deep neural network
$\mathcal{M}_\Theta$, the training strategy for finding a suitable $\mathcal{M}_\Theta$, and the method
for applying $\mathcal{M}_\Theta$ to retrieve a good approximate support matrix.

\subsubsection{Network architecture}\label{S321}

\begin{figure}[htbp]
\centering
\includegraphics[width=.85\textwidth]{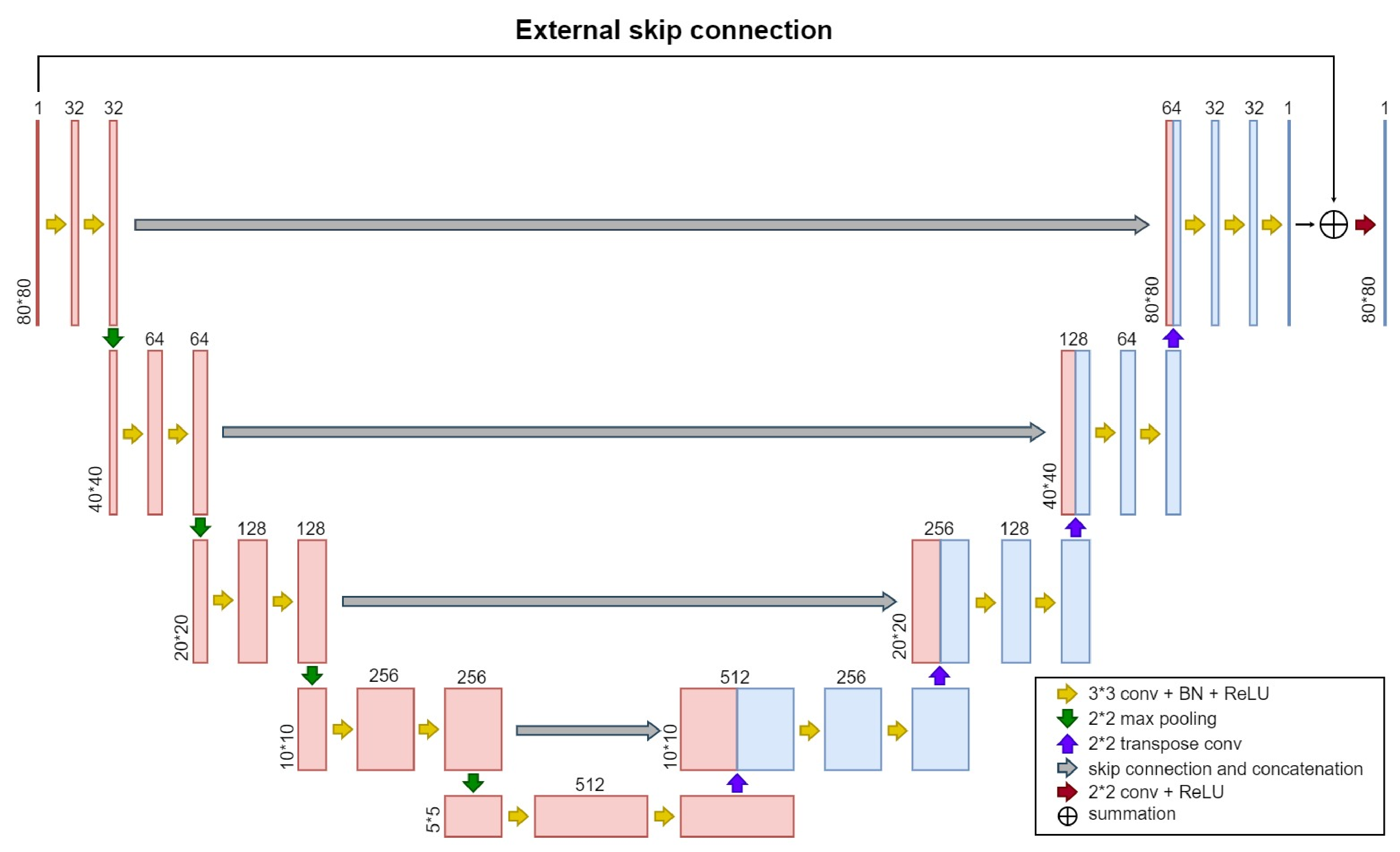}
\caption{The architecture of $\mathcal{M}_\Theta$. Each red and blue item represents a volume (also called
multichannel feature map). The number of channels is shown at the top of the volume, and the length and width
are provided at the lower-left edge of the volume. The arrows denote the different operations, which are
explained at the lower-right corner of the figure.
}\label{F2}
\end{figure}

We parameterize $\mathcal{M}_\Theta$ by a convolutional neural network called U-Net, which typically has a
U-shaped structure. The original version of U-Net was first proposed in \cite{R15} for biomedical image
segmentation. In this paper, we use a modified version of U-Net (see Figure \ref{F2}), which is similar to
but different from the one used in \cite{LZZ24}.
We note that there are only two differences between the neural network in this paper and the one used in \cite{LZZ24}.
The first difference is that the input and output of $\mathcal{M}_\Theta$ are volumes of size $N\times N\times 1$ with $N=80$ being the same as in Section \ref{S2}, which also leads to changes in the sizes of volumes of hidden layers in our network.
The second difference is that, in our network, we add a rectified linear unit (ReLU) behind the $ (2\times 2) $ convolution with $ (1\times 1) $ convolution stride to obtain the final output (see red right arrow in Figure \ref{F2}).
For more details of our neural network $\mathcal{M}_\Theta$ and the explanations of Figure \ref{F2}, the reader is referred to \cite[Section 4.1]{LZZ24}.

\subsubsection{The training strategy and application}\label{S3_2_2}

We aim to train an appropriate $\mathcal{M}_\Theta$ that maps the normalization of the imaging matrix
$\mathcal{N}(\boldsymbol{I})$ to its corresponding support matrix $\boldsymbol{S}$, where $\boldsymbol{I}$
and $\boldsymbol{S}$ denote the imaging matrix and the support matrix given as in Subsection \ref{S3_1} and
Section \ref{S2}, respectively. Here, $\mathcal{N}$ is a normalization operator defined by
$\mathcal{N}(\boldsymbol{f}):=\boldsymbol{f}/\|\boldsymbol{f}\|_{\max}$ for any
$\boldsymbol{f}:=(\boldsymbol{f}_{ij})\in\C^{N\times N}$ with the norm
$\|\boldsymbol{f}\|_{\max}:=\max_{1\leq i,j\le N}|\boldsymbol{f}_{ij}|$.
As used in \cite{LZZ24}, we hope that applying $\mathcal{N}$ to the imaging matrix $\boldsymbol{I}$
can make $\mathcal{M}_\Theta$ focus on extracting the a priori information of the shape and location
(i.e., the support) of the unknown contrast contained in the corresponding imaging matrix $\boldsymbol{I}$,
and thus leads to a good approximation to the corresponding support matrix $\boldsymbol{S}$.
During the training stage, we first generate a sample set of exact contrast matrices
$\{\boldsymbol{m}_i\}_{i=1}^T$ with $T\in\mathbb{N}^+$.
Then we simulate the noisy far-field data of $\{\boldsymbol{m}_i\}_{i=1}^T$ for the direct imaging method
(see Subsection \ref{S3_1}) to generate their corresponding imaging matrices $\{\boldsymbol{I}_i\}_{i=1}^T$.
We further compute the support matrices $\{\boldsymbol{S}_i\}_{i = 1}^T$ of the sample set
$\{\boldsymbol{m}_i\}_{i=1}^T$. We now train the deep neural network $\mathcal{M}_\Theta$ using the following
dataset with Xavier initialization \cite{GB10}:
\begin{equation}\label{27}
S:=\{(\mathcal{N}(\boldsymbol{I}_i),\boldsymbol{S}_i)\}_{i=1}^T,
\end{equation}
where $\mathcal{M}_\Theta$ is trained to map $\mathcal{N}(\boldsymbol{I}_i)$ to $\boldsymbol{S}_i$,
$i=1,2,\ldots,T$, and the loss function is given as
\begin{equation}\label{28}
\mathcal{L}(\Theta):=\sum_{i=1}^T||\mathcal{M}_\Theta(\mathcal{N}(\boldsymbol{I}_i))-\boldsymbol{S}_i||^2.
\end{equation}
Here and throughout this paper, $\|\cdot\|$ denotes the Frobenius norm of a matrix.
Finally, the desired deep neural network $\mathcal{M}_{\widehat{\Theta}}$ is obtained after training
for $t$ epochs.

The application method of our well-trained $\mathcal{M}_{\widehat{\Theta}}$ is described below.
For an unknown contrast with support matrix $\boldsymbol{S}$, we first obtain its measurement matrix
$\boldsymbol{u}^{\infty,\delta}$ (see Section \ref{S2}) for the direct imaging method and then
its corresponding imaging matrix $\boldsymbol{I}$ (see Subsection \ref{S3_1}).
Further, we compute $\boldsymbol{S}_A:=\mathcal{M}_{\widehat{\Theta}}(\mathcal{N}(\boldsymbol{I}))$
to approximate the exact support matrix $\boldsymbol{S}$.
However, $\boldsymbol{S}_A$ is not directly usable in practice since $\boldsymbol{S}_A$
is not necessarily a binary matrix, which better characterizes the approximate support area.
To address this issue, we introduce the operator $\mathcal{S}_\g$ as follows:
\begin{equation}\label{29}
\mathcal{S}_\g(\boldsymbol{f}):=(\mathcal{S}_\g(\boldsymbol{f}_{ij}))\in\R^{N\times N}\;\;\mbox{with}\;\;
\mathcal{S}_\gamma(\boldsymbol{f}_{ij})=\left\{ \begin{array}{ll}
1, & \boldsymbol{f}_{ij}>\gamma,\\
0, & \boldsymbol{f}_{ij}\leq\gamma,
\end{array} \right.
\quad \forall\boldsymbol{f}:=(\boldsymbol{f}_{ij})\in\mathbb{R}^{N\times N}.
\end{equation}
We then use $\wid{\boldsymbol{S}}:=\mathcal{S}_\g(\mathcal{M}_{\wi{\Theta}}(\mathcal{N}(\boldsymbol{I})))$
to approximate the exact support matrix $\boldsymbol{S}$, where $\g\in (0,1)$ will be chosen carefully
in the numerical reconstruction process (see Subsection \ref{S5_1_2} for the choice of $\g$).
Following the above procedure, we are able to obtain a reasonable approximation of the support of the unknown
contrast, including the high contrast case.

\section{Inversion algorithms based on the learned regularizers}\label{S4}
\setcounter{equation}{0}

In this section, we present two inversion algorithms for solving the inverse problem (IP) which incorporate
the a priori information of the shape and location (i.e., the support) of the unknown contrast as certain
regularization strategies, where the a priori information is retrieved from the direct imaging method
(see Section \ref{S3}) by our well-trained deep neural network $\mathcal{M}_{\widehat{\Theta}}$ (
see Subsection \ref{S3_2}).

\subsection{Learned projected Landweber method}\label{S4_1}

The projected Landweber method has been extensively studied for linear and nonlinear ill-posed problems
(see, e.g., \cite{KN06,PB97,N88}).
To incorporate the a priori information of the support of the unknown contrast into the projected Landweber
method, we first use the neural network $ \mathcal{M}_{\widehat{\Theta}} $ to extract an approximate support
matrix $\widetilde{\boldsymbol{S}}$ (see Subsection \ref{S3_2_2} for more details) from the direct imaging
method. Then the projected Landweber method aims to seek solutions lying on the set
$\mathcal{C}:=\{\widetilde{\boldsymbol{S}}\odot\boldsymbol{f}:\boldsymbol{f}\in\mathbb{C}^{N\times N}\}$,
where $\odot$ denotes the element-wise multiplication of two matrices of the same size.
It should be noted that the set $\mathcal{C}$ satisfies that for $\boldsymbol{f}=(\boldsymbol{f}_{ij})\in\mathcal{C}$,
\begin{equation}\label{10}
\{(i,j):\boldsymbol{f}_{ij}\neq0\}\subset\{(i,j):\widetilde{\boldsymbol{S}}_{ij} = 1\}.
%\quad \forall \boldsymbol{f}\in\mathcal{C}.
\end{equation}%
Further, for any $\boldsymbol{f}\in\C^{N\times N}$
define $\mathcal{P}_\mathcal{C}(\boldsymbol{f}):=\wid{\boldsymbol{S}}\odot\boldsymbol{f},$
which is an orthogonal projection operator from $\C^{N\times N}$ onto $\mathcal{C}$.
Then our projected Landweber iteration is given as follows:
\begin{equation}\label{30}
\begin{aligned}
\boldsymbol{m}^\delta_{i+1}
&=\mathcal{P}_\mathcal{C}\left(\boldsymbol{m}^\delta_i-\mu[\boldsymbol{F}'(\boldsymbol{m}^\delta_i)]^*
(\boldsymbol{F}(\boldsymbol{m}^\delta_i) - \boldsymbol{u}^{\infty,\delta})\right)\\
&=\widetilde{\boldsymbol{S}}\odot\left(\boldsymbol{m}^\delta_i-\mu[\boldsymbol{F}'(\boldsymbol{m}^\delta_i)]^*
(\boldsymbol{F}(\boldsymbol{m}^\delta_i) - \boldsymbol{u}^{\infty,\delta})\right),
\end{aligned}
\end{equation}
where $\boldsymbol{m}^\delta_i$ and $\boldsymbol{m}^\delta_{i+1}$ are the approximation to the unknown
contrast at the $i$-th and $(i+1)$-th iterations, respectively.
Here, the superscript $ \delta $ indicates the dependence on the noise level,
$[\boldsymbol{F}'(\boldsymbol{m}^\delta_i)]^*$ is the adjoint of $\boldsymbol{F}'(\boldsymbol{m}^\delta_i)$
and $\mu>0$ is the stepsize.
Since the projection operator $\mathcal{P}_\mathcal{C}$ can project an approximate contrast matrix into $\mathcal{C}$,
then, by \eqref{10} it is reasonable to expect that
$\mathcal{P}_\mathcal{C}$ can be helpful for choosing a desirable solution when the approximate support
matrix $\widetilde{\boldsymbol{S}}$ is a good approximation to the exact support matrix $\boldsymbol{S}$,
which is the training purpose of our deep neural network $\mathcal{M}_{\widehat{\Theta}}$
(see Subsection \ref{S3_2}).
Since the a priori information incorporated in the projected Landweber method is learned by our deep neural network $\mathcal{M}_{\widehat{\Theta}}$, we also call it the learned projected Landweber method.

We now describe our learned projected Landweber method for the problem (IP).
In what follows, for any wave number $k$, we rewrite the direct operator $\boldsymbol{F}$ and the noisy
far-field data $\boldsymbol{u}^{\infty, \delta}$ in \eqref{5} as $\boldsymbol{F}_k$ and
$\boldsymbol{u}^{\infty, \delta}_k $, respectively, to indicate the dependence on the wave number $k$.
We first use the noisy far-field data $\boldsymbol{u}_{k_0}^{\infty,\delta}$ measured at the wave number $k_0$
to obtain its imaging matrix $\boldsymbol{I}$ (see Subsection \ref{S3_1}) by the direct imaging method,
and then compute the approximate support matrix
$\wid{\boldsymbol{S}}:= \mathcal{S}_\gamma(\mathcal{M}_{\wi{\Theta}}(\mathcal{N}(\boldsymbol{I})))$
(see Subsection \ref{S3_2_2}). For the learned projected Landweber method, we make use of the noisy far-field
data $\boldsymbol{u}_k^{\infty,\delta}$ measured at the wave number $k$.
In the reconstruction process, we choose $\mu$ as the stepsize and set the initial guess of $\boldsymbol{m}$
to be $0$. Let $N_p$ be the total iteration number.
Then the learned projected Landweber method is presented in Algorithm \ref{P} for the problem (IP).
For simplicity, we call Algorithm \ref{P} as \textit{Learned Projected Algorithm} in the remaining part
of the paper. See Section \ref{S5} for the performance of this algorithm.

\begin{algorithm}[htbp]
\caption{Learned projected Landweber method}\label{P}

\textbf{Input: }$\boldsymbol{F}_k$, $\boldsymbol{u}_k^{\infty,\delta}$, $\boldsymbol{u}_{k_0}^{\infty,\delta}$,
$\gamma$, $\mathcal{M}_{\widehat{\Theta}}$, $\mu$, $N_p$

\textbf{Output:} the final approximate contrast for $\boldsymbol{m}$

\textbf{Initialize:} $i = 0$, $\boldsymbol{m}_0^\delta = 0$

\begin{algorithmic}[1]
\item Compute the imaging matrix $\boldsymbol{I}$ with the far-field data $\boldsymbol{u}_{k_0}^{\infty,\delta}$
\item Compute $\wid{\boldsymbol{S}}:=\mathcal{S}_\g(\mathcal{M}_{\wi{\Theta}}(\mathcal{N}(\boldsymbol{I})))$
\item \textbf{while} $i<N_p$ \textbf{do}
\item \begin{itemize}
\item[]  $\boldsymbol{m}_{i+1}^\delta =\widetilde{\boldsymbol{S}}\odot\left(\boldsymbol{m}^\delta_i
         - \mu[\boldsymbol{F}_k'(\boldsymbol{m}^\delta_i)]^*(\boldsymbol{F}_k(\boldsymbol{m}^\delta_i)
         - \boldsymbol{u}^{\infty,\delta}_k)\right)$
\end{itemize}
\item \begin{itemize}
\item[]  $ i\gets i+1 $	
\end{itemize}
\item \textbf{end while}
\item Set the final approximate contrast to be $\boldsymbol{m}_{N_p}^\delta$.
\end{algorithmic}
\end{algorithm}

\begin{remark}\label{R1} {\rm
If the projection operator $\mathcal{P_{\mathcal{C}}}$ is replaced by the identity map, then the projected Landweber
iteration \eqref{30} reduces to the standard Landweber iteration \cite{L51}, which is called
\textit{Landweber Algorithm} in this paper. To demonstrate the benefits provided by our deep neural network
$\mathcal{M}_{\widehat{\Theta}}$, we will compare this algorithm with the proposed reconstruction algorithms
in Subsection \ref{S5_2}.
}
\end{remark}

\subsection{Learned variational regularization method}\label{S4_2}

One of the most popular approaches for solving the inverse problem (IP) is the variational regularization method.
This approach reformulates the inverse problem (IP) as the following optimization problem
\begin{equation}\label{18}
\underset{\boldsymbol{m}\in \mathbb{C}^{N\times N}}{\arg\min}\dfrac{1}{2}\left|\left|\boldsymbol{F}(\boldsymbol{m})
- \boldsymbol{u}^{\infty,\delta}\right|\right|^2 + \lambda\mathcal{R}(\boldsymbol{m}),
\end{equation}
where $\mathcal{R}:\mathbb{C}^{N\times N}\to [0,+\infty)$ is an appropriate regularization functional which encodes
certain a priori information of the exact contrast matrix and should be given in advance, and
%penalizes unfeasible solutions, and
$\lambda>0$ is a regularization parameter that governs the influence of the a priori knowledge encoded
by the regularization functional on the need to fit data and also needs to be provided in advance.
In this section, we will design the regularization functional $\mathcal{R}$ by using the a priori information of
the support of the unknown contrast. Precisely, we first use the deep neural network $\mathcal{M}_{\widehat{\Theta}}$
to extract the approximate support matrix $\widetilde{\boldsymbol{S}}$ (see Subsection \ref{S3_2_2} for more
details) from the direct imaging method, and then define the regularization functional $\mathcal{R}$ as
$$
\mathcal{R}(\boldsymbol{f}):=(1/2)||\boldsymbol{f}-\widetilde{\boldsymbol{S}}\odot\boldsymbol{f}||^2,\quad
\forall \boldsymbol{f}\in\mathbb{C}^{N\times N}.
$$
%in \eqref{18}.
Note that the regularization functional $\mathcal{R}$ defined above is expected to penalize large values
of an approximate contrast matrix outside the index set $\{(i,j): \widetilde{\boldsymbol{S}}_{ij}=1\}$ and
thus hopefully can lead to good reconstruction results if $\widetilde{\boldsymbol{S}}$ is a good
approximation to the exact support matrix $\boldsymbol{S}$, which is actually the training purpose of
our deep neural network $\mathcal{M}_{\widehat{\Theta}}$ (see Subsection \ref{S3_2}).
In order to solve the optimization problem \eqref{18} with $\mathcal{R}$ defined above, we use the gradient
descent method (GD), which updates at each iteration as follows (see \cite[Section 4.5]{C19}):
\begin{equation}\label{20}
\boldsymbol{m}^\delta_{i+1} = \boldsymbol{m}^\delta_i
- \mu\left([\boldsymbol{F}'(\boldsymbol{m}^\delta_i)]^*(\boldsymbol{F}(\boldsymbol{m}^\delta_i)
- \boldsymbol{u}^{\infty,\delta}) + \lambda(\boldsymbol{m}_i^\delta
- \boldsymbol{m}_i^\delta\odot\widetilde{\boldsymbol{S}})\right),
\end{equation}
where $\boldsymbol{m}^\delta_i$ and $\boldsymbol{m}^\delta_{i+1}$ are the approximation to the unknown
contrast at the $i$-th and $(i+1)$-th iterations, respectively, with the superscript $\delta$ indicating
the dependence on the noise level, and $\mu>0$ is the stepsize of GD.
Since the a priori information encoded in the regularization functional $\mathcal{R}$ is learned by the
deep neural network $\mathcal{M}_{\widehat{\Theta}}$, the iteration method \eqref{20} is then called
the learned variational regularization method.

We now describe the above learned variational regularization method for the problem (IP).
We use the same notations and measured data as in Section \ref{S4_1}.
In the reconstruction process, we choose $\mu$ as the stepsize and set the initial guess of
$\boldsymbol{m}$ to be $0$. Let $N_v$ be the total iteration number.
Then the learned variational regularization method is presented in Algorithm \ref{V} for the problem (IP).
For simplicity, we call Algorithm \ref{V} as \textit{Learned Variational Algorithm} in the remaining part
of the paper. See Section \ref{S5} for the performance of this algorithm.

\begin{algorithm}[htbp]
\caption{Learned variational regularization method}\label{V}

\textbf{Input:} $\boldsymbol{F}_k$, $\boldsymbol{u}_k^{\infty,\delta}$, $\boldsymbol{u}_{k_0}^{\infty,\delta}$,
$\gamma$, $\mathcal{M}_{\widehat{\Theta}}$, $\mu$, $\lambda$, $N_v$

\textbf{Output:} the final approximate contrast for $\boldsymbol{m}$

\textbf{Initialize:} $ i = 0 $, $ \boldsymbol{m}_0^\delta = 0 $

\begin{algorithmic}[1]
\item Compute the imaging matrix $\boldsymbol{I}$ with the far-field data $\boldsymbol{u}_{k_0}^{\infty,\delta}$
\item Compute $\wid{\boldsymbol{S}}:=\mathcal{S}_\gamma(\mathcal{M}_{\wi{\Theta}}(\mathcal{N}(\boldsymbol{I})))$
\item \textbf{while} $i<N_v$ \textbf{do}
\item \begin{itemize}
\item[]  $\boldsymbol{m}^\delta_{i+1} = \boldsymbol{m}^\delta_i - \mu\left([\boldsymbol{F}_k'(\boldsymbol{m}^\delta_i)]^*(\boldsymbol{F}_k(\boldsymbol{m}^\delta_i) - \boldsymbol{u}^{\infty,\delta}_k) + \lambda(\boldsymbol{m}_i^\delta - \boldsymbol{m}_i^\delta\odot\widetilde{\boldsymbol{S}})\right)$
\end{itemize}
\item \begin{itemize}
\item[]  $ i\gets i+1 $	
\end{itemize}
\item \textbf{end while}
\item Set the final approximate contrast to be $\boldsymbol{m}_{N_v}^\delta$.
\end{algorithmic}
\end{algorithm}

\begin{remark}\label{R3}\rm
While the deep neural network is used to learn regularizers for our inversion algorithms, there is another popular method that applies a deep neural network as a post-processing step to solve inverse problems.
This kind of method is known as the post-processing method.
Roughly speaking, a post-processing method first applies a traditional inversion method to generate an initial reconstruction, which is then refined by passing it through a trained deep neural network. For a detailed discussion of post-processing methods in various inverse problems, we refer to \cite[Section 5.1.5]{AMOS19}.
In Subsection \ref{S5_2}, we will give some numerical results for the inverse problem (IP) by a post-processing method.
Specifically, given any exact contrast matrix $\boldsymbol{m}$, we first compute its initial approximation by using \textit{Landweber Algorithm} (see Remark \ref{R1}), then refine this approximation through a trained neural network, which is trained to map the approximations generated by \textit{Landweber Algorithm} to their corresponding exact contrast matrices.
For fair comparison, we also use $\mathcal{M}_\Theta$ given in Subsection \ref{S321} for the trained neural network in the above post-processing method; see Subsection \ref{S5_1_2} for the training strategy of the neural network.
For simplicity, we call this algorithm as \textit{Post-processing Algorithm} in this paper.
The comparisons between our inversion algorithms (i.e., \textit{Learned Projected Algorithm} and \textit{Learned Variational Algorithm}) and \textit{Post-processing Algorithm} will be given in Subsection \ref{S5_2}.
\end{remark}

\section{Numerical experiments}\label{S5}
\setcounter{equation}{0}

In this section, we present numerical examples to demonstrate the effectiveness of the proposed inversion
algorithms based on the deep learning method (i.e., \textit{Learned Projected Algorithm} and
\textit{Learned Variational Algorithm}) for the problem (IP).
The experimental setup is given in Subsection \ref{S5_1} for numerical experiments.
The performance of \textit{Learned Projected Algorithm} and \textit{Learned Variational Algorithm} are shown
in Subsection \ref{S5_2}. In order to evaluate the role played by the a priori information of the shape and
location (i.e., the support) of the unknown contrast encoded in our algorithms, we will also compare these
two algorithms with \textit{Landweber Algorithm} (see Remark \ref{R1}) and \textit{Post-processing Algorithm} (see Remark \ref{R3}) in Subsection \ref{S5_2}.
To show the robustness of the proposed algorithms with respect to noise, they are tested with different
noise levels in Subsection \ref{S5_3}.

\subsection{Experimental setup}\label{S5_1}

The training process is performed on COLAB (Tesla P100 GPU, Linux operating system) and is implemented on PyTorch,
while the computations of direct scattering problem, \textit{Learned Projected Algorithm} and
\textit{Learned Variational Algorithm} are implemented by Python 3.9 on a desktop computer
(Intel Core i7-10700 CPU (2.90 GHz), 32 GB of RAM, Ubuntu 20.04 LTS).

\subsubsection{Simulation setup for the scattering problem}\label{S5_1_1}

As mentioned in Section \ref{S2}, the support of the unknown contrast is assumed to lie in a disk
$B_\rho\subset C_\rho$ with $\rho>0$. Without loss of generality, we choose $\rho=3$.
For the measured far-field data $\boldsymbol{u}^{\infty,\delta}_k$ (see Subsections \ref{S4_1} and \ref{S4_2}),
which is used in the iterative process, we set the number of incident directions $Q_1=16$ and the number
of measured directions $P_1= 32$.
For the measured far-field data $\boldsymbol{u}_{k_0}^{\infty,\delta}$ used in the direct imaging method,
we set the number of incident directions $ Q_2 = 64 $ and the number of measured directions $ P_2 = 128 $.
To generate the synthetic far-field data, we use the method discussed in Remark \ref{R2} with $N=320$.
The noisy far-field data $\boldsymbol{u}^{\infty,\delta}(\hat{x}_p,d_q)$, $p=1,\ldots,P_i,\;q=1,\ldots,Q_i$
($ i = 1,2 $), are given as mentioned in Section \ref{S2}. In the training stage, we choose $\delta = 5\%$
for both $\boldsymbol{u}^{\infty,\delta}_k $ and $\boldsymbol{u}_{k_0}^{\infty,\delta}$.

\subsubsection{Parameter setting for inversion algorithms}\label{S5_1_2}

For the parameters in \textit{Learned Projected Algorithm}, we choose $N=80$, $N_p=100$, wave numbers $k=2$
and $k_0 = 15$, stepsize $\mu = 0.5$ and the threshold $\gamma = 0.1$.
For \textit{Learned Variational Algorithm}, we choose the same parameters as those in \textit{Learned Projected Algorithm}, and we additionally choose $ N_v = 100 $ and the regularization parameters $ \lambda = 1 $.
Further, in order to obtain the well-trained deep neural network $\mathcal{M}_{\widehat{\Theta}}$ for these two algorithms, we train $\mathcal{M}_\Theta$ by minimizing the loss function $\mathcal{L}$ with the epochs $t=20$,
using the Adam optimizer \cite{KB14} with batch size $10$ and learning rate $10^{-3}$.
For \textit{Landweber Algorithm} given in Remark \ref{R1}, the total iteration number and stepsize are chosen to be the same as those in \textit{Learned Projected Algorithm}.
For \textit{Post-processing Algorithm} given in Remark \ref{R3}, we train the deep neural network by minimizing the standard mean squared error loss (MSE) with $ 200 $ epochs, using Adam optimizer \cite{KB14} with batch size $ 10 $ and learning rate $ 10^{-3} $.

\subsubsection{Data generation for the neural network}\label{S5_1_3}

In this paper, we consider the exact contrast $m$ consisting of several disjoint inhomogeneous media with
the shape of their supports being ellipse.
Precisely, let $U(a,b)$ be uniform distribution in $ [a,b] $ and
\begin{equation}\label{21}
R_\theta := \left(\begin{array}{cc}
	\cos\theta & -\sin\theta \\
	\sin\theta & \cos\theta
\end{array}\right)	
\end{equation}
be the rotation transformation, we randomly select two or three disjoint ellipses
$e:=\{(x_0,y_0)+R_{\theta}(x,y):(x/a)^2+(y/b)^2\leq 1, a>b\}\subset B_\rho$, where $(x_0,y_0)$ is
the center of $e$, $ \theta $ is the rotation angle of $ e $ which is sampled from $ U(0,2\pi) $, and
the semi-major axis $ a $ and semi-minor axis $ b $ are sampled from $ U(0.6,1.2) $ and $ U(0.3,0.6) $,
respectively. Then we consider the dataset consisting of the exact contrast $ m $ supported in these ellipses
with the definition $ m(x) = c,\; x\in e $, where $ c $ are sampled from $ U(1,3) $ for each ellipse.
This dataset will be used to train the proposed deep neural network $\mathcal{M}_\Theta$
(see Subsection \ref{S5_2}) and called \textit{Ellipse Dataset} in the present paper.

\subsubsection{Evaluation criterion for inversion algorithms}

In order to quantitatively evaluate the reconstruction performance of an inversion algorithm for
the problem (IP), we introduce an error function to measure the difference between the exact refractive
index and the approximate refractive index generated by an inversion algorithm.
As mentioned before, for the exact contrast $m(x)$, $\boldsymbol{m}=(\boldsymbol{m}_{ij})\in\mathbb{C}^{N\times N}$
is the exact contrast matrix with $\boldsymbol{m}_{ij} = m(x_{ij})$ and
$\widehat{\boldsymbol{m}} = (\widehat{\boldsymbol{m}}_{ij}) $ is the output of an inversion algorithm,
which is the approximation of $\boldsymbol{m}$. Here, $x_{ij}$ ($i,j=1,2,\ldots,N$) are the points introduced
at the end of Section \ref{S2}. Accordingly, $\boldsymbol{n}=(\boldsymbol{n}_{ij}):=\boldsymbol{m}+1$ is
the discretization of the refractive index $n(x)=m(x)+1$ and $\wi{\boldsymbol{n}}=(\wi{\boldsymbol{n}}_{ij})
:= \wi{\boldsymbol{m}} + 1 $ is the approximation of $ \boldsymbol{n} $.
We now define the relative error function $\mathcal{E}$ between $\boldsymbol{n}$ and $\widehat{\boldsymbol{n}}$
as follows:
\begin{displaymath}
\mathcal{E}(\boldsymbol{n},\widehat{\boldsymbol{n}}):= \dfrac{\|\boldsymbol{n} - \widehat{\boldsymbol{n}}\|}{\|\boldsymbol{n}\|},
\end{displaymath}
where $\|\cdot\|$ denotes the Frobenius norm of a matrix, as mentioned in Subsection \ref{S3_2_2}.

\subsection{Performance of the proposed inversion algorithms}\label{S5_2}

We use \textit{Ellipse Dataset} (see Subsection \ref{S5_1_3}) to train our deep neural network
$\mathcal{M}_\Theta$. During the training process, $1800$ samples are used for training $\mathcal{M}_\Theta$
with the training strategy mentioned in Subsection \ref{S3_2_2} and $200$ samples are used to validate the training performance.
For fairness, we adopt the same training samples for the deep neural network used in \textit{Post-processing Algorithm}.

We first investigate the influence of the value of different contrasts on \textit{Landweber Algorithm}, \textit{Post-processing Algorithm},
\textit{Learned Projected Algorithm} and \textit{Learned Variational Algorithm} quantitatively.
To do this, we consider three cases with the noise level $\delta=5\%$, which are denoted as Cases 1.1, 1.2 and 1.3.
In all three cases, we generate $100$ samples from the \textit{Ellipse Dataset} to represent the exact
contrast matrices.
For Cases 1.1, 1.2 and 1.3, we set the norm $\|\cdot\|_{\max}$ of each exact contrast matrix to be $2,3$
and $4,$ respectively.
Table \ref{T1} presents the average values of the relative error $\mathcal{E}$ for the four algorithms
mentioned above for these three cases.
Figure \ref{F1} presents the reconstruction results of several samples from the three cases.
Each row of Figure \ref{F1} shows the normalization of the imaging matrix (see step 1 of Algorithms \ref{P} and \ref{V}),
the approximate support generated by $\mathcal{M}_{\widehat{\Theta}}$ (see step 2 of Algorithms \ref{P} and \ref{V}), the numerical reconstructions by the
four algorithms and the ground truth for one sample. It can be seen in Table \ref{T1} and Figure \ref{F1}
that the proposed inversion algorithms outperform the \textit{Landweber Algorithm} and \textit{Post-processing Algorithm}, especially in the high contrast case, which shows the advantage
offered by our deep neural network $\mathcal{M}_{\widehat{\Theta}}$.
Moreover, our experiments show that the proposed inversion algorithms give satisfactory results even for the case when the value of the contrasts for testing is higher than those for training (see the last column of Table \ref{T1} as well as the fifth and sixth rows of Figure \ref{F1}).
This may be because our deep neural network $\mathcal{M}_{\widehat{\Theta}}$ can extract
the a priori information of the support from the direct imaging method for various values of different contrasts, which is critical a priori information for good reconstruction.

\begin{figure}[htbp]
\centering
\includegraphics[width=1.\textwidth]{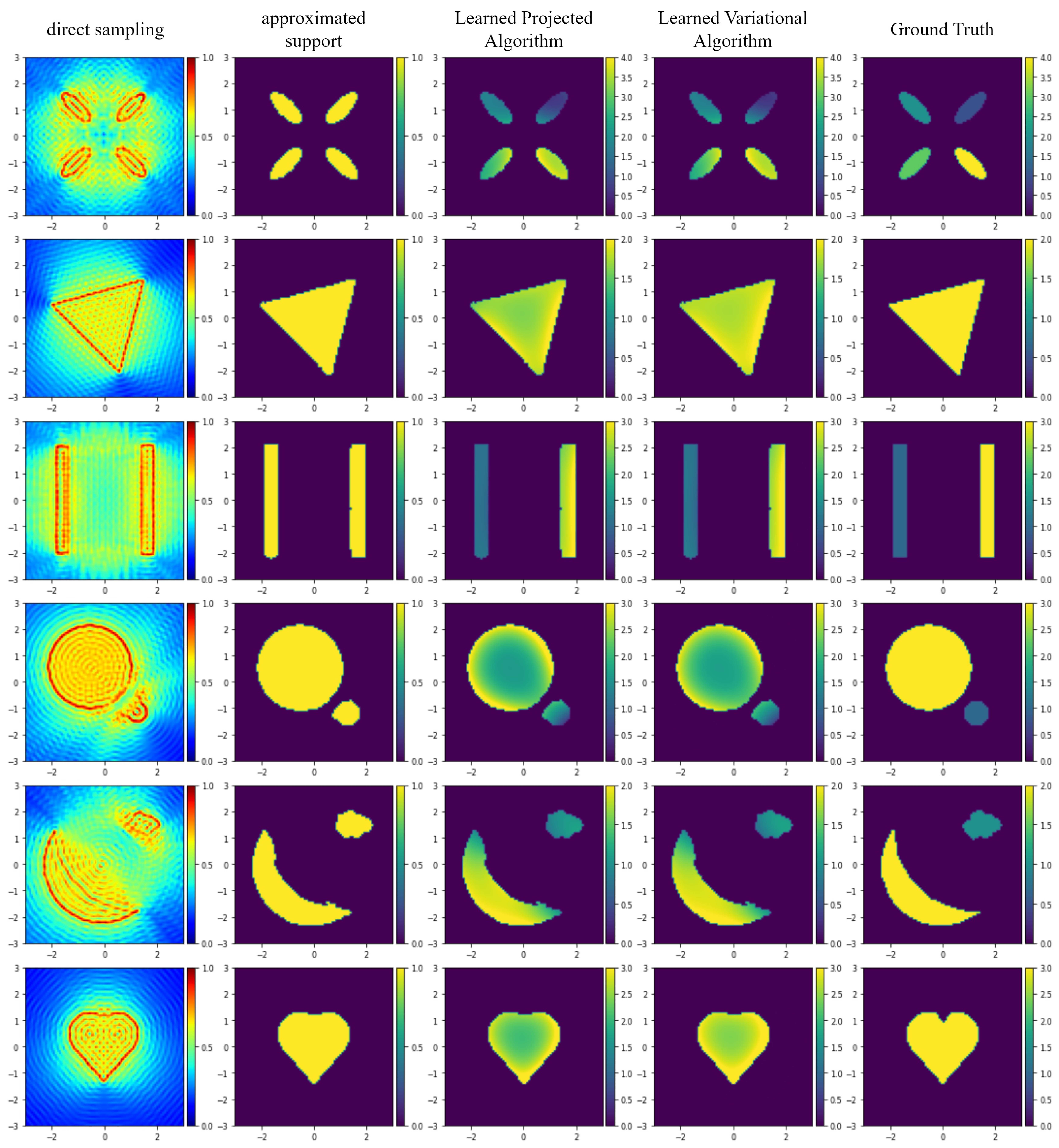}
\caption{Reconstruction results by \textit{Landweber Algorithm}, \textit{Post-processing Algorithm}, \textit{Learned Projected Algorithm}
and \textit{Learned Variational Algorithm} for the test samples outside
\textit{Ellipse Dataset}.
Each row presents the normalization of the imaging matrix, the approximate support generated by $\mathcal{M}_{\widehat{\Theta}}$, the reconstructions by \textit{Landweber Algorithm}, \textit{Post-processing Algorithm},
\textit{Learned Projected Algorithm}, \textit{Learned Variational Algorithm} and the ground truth for one sample.}\label{F3}
\end{figure}

We now test the generalization ability of our inversion algorithms with the exact contrast matrix out
of \textit{Ellipse Dataset}, where the noise level is set to be $ \delta = 5\% $.
The reconstruction results by \textit{Landweber Algorithm}, \textit{Post-processing Algorithm}, \textit{Learned Projected Algorithm} and \textit{Learned Variational Algorithm}
are shown in Figure \ref{F3}, where each row presents the normalization of the imaging matrix (see step 1 of Algorithms \ref{P} and \ref{V}), the approximate support generated by $\mathcal{M}_{\widehat{\Theta}}$ (see step 2 of Algorithms \ref{P} and \ref{V}), the reconstructions by the four algorithms, and
the ground truth for one sample. The reconstruction results in Figure \ref{F3} demonstrate
the good generalization ability of the proposed algorithms, which can recover unknown contrasts with various
support shapes and contrast values, even though they are not sampled from the \textit{Ellipse Dataset}.
Regarding this, we believe that our training strategy makes $\mathcal{M}_{\widehat{\Theta}}$ focus on learning the a priori information of the support of the unknown contrast, which leads to a satisfactory generalization ability of  $\mathcal{M}_{\widehat{\Theta}}$ for various support shapes.
It should be also noted that, although our deep neural network $\mathcal{M}_{\widehat{\Theta}}$ is trained only on samples with piecewise constant contrasts (i.e., \textit{Ellipse Dataset}), our inversion algorithms can provide satisfactory reconstructions for test samples with piecewise smooth contrasts (see the third, forth and fifth rows in Figure \ref{F3}).
Remarkably, it is shown in Figure \ref{F3} that for some test samples, although the approximate supports computed by $\mathcal{M}_{\widehat{\Theta}}$ are not very accurate, our inversion algorithms can still
achieve stable reconstructions.

\begin{table}[htpb]
	\centering
	\begin{tabular}{|p{5.5cm}<{\centering}|p{2cm}<{\centering}|p{2cm}<{\centering}|p{2cm}<{\centering}|}
		\hline
		& Case 1.1 & Case 1.2 & Case 1.3\\
		\hline
		\textit{Landweber Algorithm} & 24.7\% & 35.8\% & 45.1\% \\
		\hline
		\textit{Post-processing Algorithm} & 10.2\% & 14.3\% & 23.2\% \\
		\hline
		\textit{Learned Projected Algorithm} & 5.7\% & 9.7\% & 15.6\% \\
		\hline
		\textit{Learned Variational Algorithm} & 5.4\% & 10.3\% & 17.0\%\\
		\hline
	\end{tabular}
\caption{The average values of the relative error $\mathcal{E}$ for the outputs of
\textit{Landweber Algorithm},  \textit{Post-processing Algorithm}, \textit{Learned Projected Algorithm} and
\textit{Learned Variational Algorithm} on the \textit{Ellipse Dataset}.
}\label{T1}
\end{table}

\begin{figure}[htbp]
\centering
\includegraphics[width=1.\textwidth]{image/performance/performance.eps}
\caption{Reconstruction results by \textit{Landweber Algorithm}, \textit{Post-processing Algorithm}, \textit{Learned Projected Algorithm}
and \textit{Learned Variational Algorithm} for the test samples from
\textit{Ellipse Dataset}. Each row presents the normalization of the imaging matrix, the approximate
support generated by $\mathcal{M}_{\widehat{\Theta}}$, the reconstructions by \textit{Landweber Algorithm}, \textit{Post-processing Algorithm},
\textit{Learned Projected Algorithm}, \textit{Learned Variational Algorithm} and the ground truth for one sample.}\label{F1}
\end{figure}

\subsection{Sensitivity to noise}\label{S5_3}

\begin{figure}[htbp]
	\centering
	\includegraphics[width=1.\textwidth]{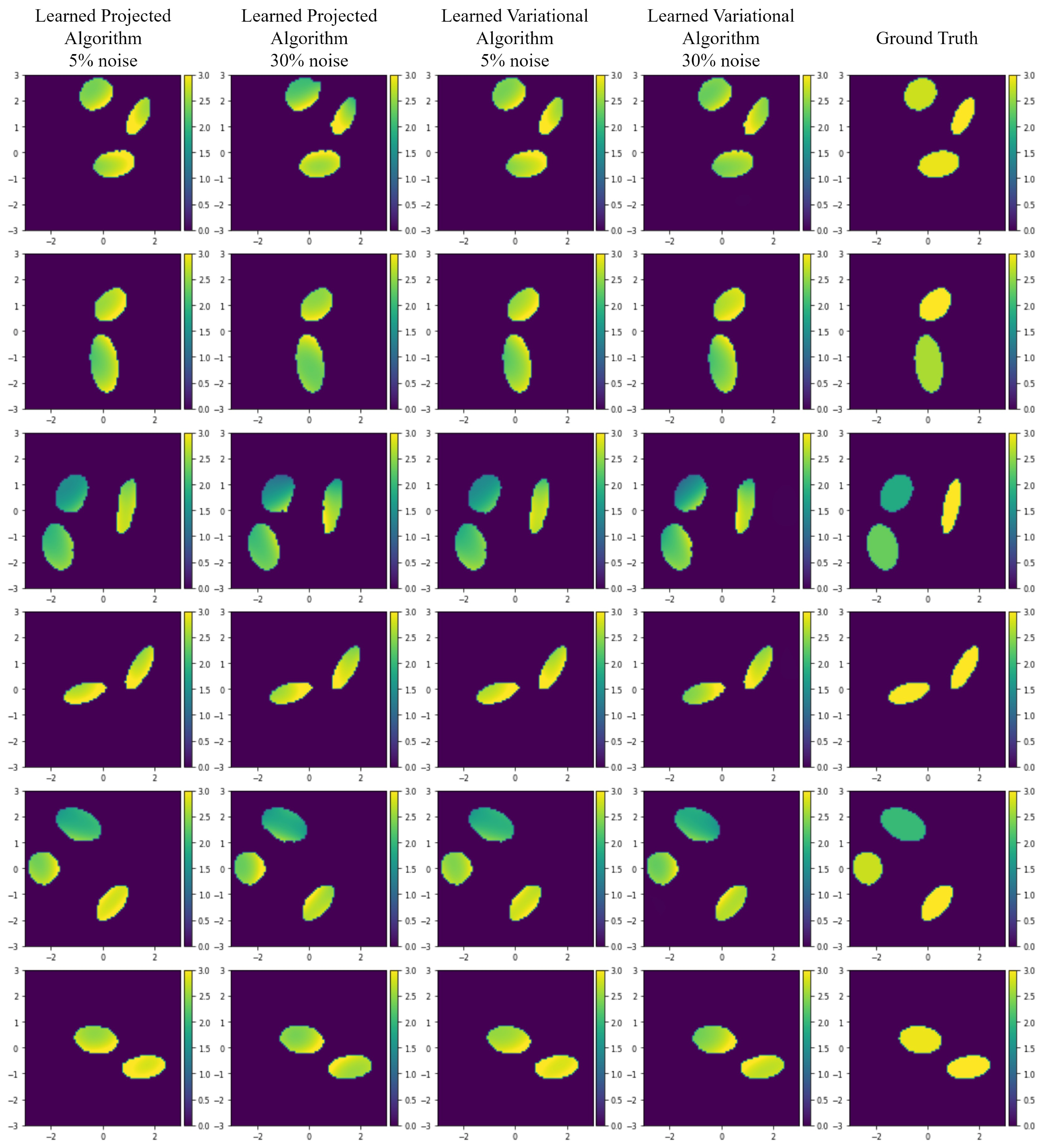}
\caption{Reconstruction results by \textit{Learned Projected Algorithm} and \textit{Learned Variational Algorithm}
with different noise levels, where the test samples are from \textit{Ellipse Dataset}.
Each row presents the reconstruction results and the ground truth for one sample.}\label{F4}
\end{figure}

In order to test the robustness of our algorithms with respect to noise, we evaluate
\textit{Learned Projected Algorithm} and \textit{Learned Variational Algorithm} under different noise settings,
where the parameters of the two algorithms are the same as in Subsection \ref{S5_1_2}. To do this, we consider
two cases with different noise levels $\delta=5\%$ and $\delta=30\%$, which are denoted as Cases 2.1 and 2.2, respectively. In both cases, we choose the same $100$ samples as in Cases 1.1, 1.2 and 1.3, to represent
the exact contrast matrices with the norm $\|\cdot\|_{\max}$ being $3.$ For these two cases, we present the
average values of the relative error $\mathcal{E}$ for the output of the two algorithms in Table \ref{T2}.
Figure \ref{F4} presents the reconstruction results of the two algorithms and the ground truth for one sample.
The reconstruction results in Table \ref{T2} and Figure \ref{F4} show that the performance of the proposed
algorithms does not degrade significantly as the noise level increases, which confirms the robustness of our
algorithms with respect to noise. Regarding this, we think this may be because the direct imaging method
is very robust to noise, which can help our deep neural network successfully extract the shape and location
information of the unknown contrast even in the presence of heavy noise.

\begin{table}[htpb]
	\centering
	\begin{tabular}{|p{6cm}<{\centering}|p{3.5cm}<{\centering}|p{3.5cm}<{\centering}|}
		\hline
		& Case 2.1 & Case 2.2\\
		\hline
		\textit{Learned Projected Algorithm} & 9.7\% & 10.5\%\\
		\hline
		\textit{Learned Variational Algorithm} & 10.3\% & 11.4\%\\
		\hline
	\end{tabular}
\caption{The average values of the relative error $\mathcal{E}$ for the outputs of \textit{Learned Projected Algorithm}
and \textit{Learned Variational Algorithm} on \textit{Ellipse Dataset}, where the noise levels for Cases 2.1
and 2.2 are set to be $\delta=5\%$ and $\delta=30\%$, respectively.}\label{T2}
\end{table}

\section{Conclusion}\label{S6}
\setcounter{equation}{0}

In this paper, we considered the inverse problem of scattering of time-harmonic acoustic waves from
inhomogeneous media in two dimensions, that is, the inverse problem of reconstructing the inhomogeneous
medium (or its contrast) from the far-field measurement data.
%including the high contrast case.
This inverse medium scattering problem is highly nonlinear and severely ill-posed, making it essential to employ regularization strategies incorporating specific a priori information of the unknown scatterer.
%However, the design of an appropriate regularization strategy relies heavily on certain a priori information
%of the unknown solution of the inverse problem which is difficult to determine.
To this end, we proposed a deep neural network $\mathcal{M}_{\widehat{\Theta}}$ to retrieve the
a priori information of the shape and location (i.e., the support) of the unknown contrast from a direct
imaging method. By incorporating the retrieved a priori information of the support of the unknown contrast
as a regularization strategy, we then proposed two inversion algorithms, \textit{Learned Projected Algorithm}
(see Subsection \ref{S4_1}) and \textit{Learned Variational Algorithm} (see Subsection \ref{S4_2})
to solve the inverse medium scattering problem.
With the aid of the a priori information of the unknown contrast provided by the deep neural
network $\mathcal{M}_{\widehat{\Theta}}$,
our inversion algorithms achieve stable reconstruction results even for the high contrast case.

Extensive numerical experiments are conducted and showed that \textit{Learned Projected Algorithm} and
\textit{Learned Variational Algorithm} perform well for the inverse scattering problem.
First, the reconstructions provided by \textit{Landweber Algorithm}, \textit{Post-processing Algorithm} and our two inversion algorithms, which were presented in Figure \ref{F1}, showed that the a priori information of the supports of the unknown contrasts retrieved by $\mathcal{M}_{\widehat{\Theta}}$ from the direct imaging method plays a vital role in our two inversion algorithms.
Second, the reconstruction results in Figure \ref{F3} showed that
the proposed two inversion algorithms have a good performance on samples outside \textit{Ellipse Dataset}
on which our deep neural network $\mathcal{M}_{\widehat{\Theta}}$ was trained,
demonstrating the satisfactory generalization ability of the proposed inversion algorithms.
In particular, it is indeed shown in Figure \ref{F3} that our two inversion algorithms can provide satisfactory reconstructions for samples with piecewise smooth contrasts, although all training samples from \textit{Ellipse Dataset} are piecewise constant.
Third, as seen from the reconstruction results in Figure \ref{F4},
the proposed two algorithms do not have a significant decrease with the added heavy noise,
indicating the robustness of the proposed inversion algorithms with respect to noise.
Regarding this, we believe this may be because the direct imaging method is strongly robust with
respect to noise and thus can help our deep neural network $\mathcal{M}_{\widehat{\Theta}}$ successfully
extract the a priori information of the support of the unknown contrast even in the presence of heavy noise.
However, it is observed in the numerical experiments that the reconstruction results of the proposed inversion
algorithms become worse when the contrast value is very large. One of the reasons may be due to
the fact that for the case when the contrast value is very large the numerical solution of
the corresponding scattering problem at each iteration may not be accurate enough, which will deteriorate
the final reconstruction results.
It is also interesting to extend our method to other challenging inverse problems,
which will be considered as a future work.

\section*{Acknowledgments}

This work was partly supported by the National Key R\&D Program of China (2024YFA1012300, 2024YFA1012303),
Beijing Natural Science Foundation (Z210001), the NNSF of China (12431016, 12271515),
and Youth Innovation Promotion Association of CAS.

\end{document}